\documentclass{amsart}
\usepackage{amsaddr}
\usepackage{hyperref}
\usepackage{etoolbox}

\usepackage{graphicx}
\usepackage{tikz}
\usepackage{tikz-cd}
\usepackage{verbatim}
\usepackage{setspace}
\usepackage{mathrsfs}

\newtheorem{theorem}{Theorem}[section]
\newtheorem{lemma}[theorem]{Lemma}

\theoremstyle{definition}

\theoremstyle{remark}

\theoremstyle{proposition}
\newtheorem{proposition}[theorem]{Proposition}

\theoremstyle{corollary}
\newtheorem{corollary}[theorem]{Corollary}

\numberwithin{equation}{section}

\DeclareMathOperator{\E}{E}
\newcommand{\abs}[1]{\lvert#1\rvert}

\newcommand*\xbar[1]{%
\hbox{%
\vbox{%
\hrule height 0.5pt 
\kern0.5ex
\hbox{%
\kern-0.1em
\ensuremath{#1}%
\kern-0.1em
}%
}%
}%
}
\newcommand\restr[2]{{
\left.\kern-\nulldelimiterspace 
#1 
\right|_{#2} 
}}

\begin{document}
\title{Preenveloping Classes of Acts}

\author{Mohanad Farhan Hamid}
\address{Department of Production and Metallurgy Engineering,
\\University of Technology-Iraq, Baghdad, Iraq
\\E-mails: {70261@uotechnology.edu.iq
\newline mohanadfhamid@yahoo.com}}

 \curraddr{}

\thanks{}





\maketitle

\section*{\bf{Abstract}}
\noindent Preenvelopes of acts over a monoid are defined by analogy with Enochs' definition of preenvelopes of modules. Provided that it is closed for pure subacts, a class of acts is shown to be preenveloping precisely when it is closed under direct products. Examples of such classes include absolutely pure, weakly f-injective and weakly p-injective acts.

\bigskip
\noindent \textbf{\keywordsname.} Pure subact, preenvelope, absolutely pure act, weakly f-injective act, weakly p-injective act.

\section{\bf{Introduction}}
\noindent Acts and $S$-acts are always right unitary over a monoid $S$. Consider systems of equations over an act $A$ (with constants from $A$). Such equations are of one of the following forms
$$xr=xs, \thinspace \thinspace xr=ys, \thinspace \thinspace \text{or} \thinspace \thinspace xr=a$$
with $r$ and $s$ $\in S$, $a$ $\in A$ and indefinite $x$ and $y$. Let $A$ be a subact of an $S$-act $B$. If every finite such system of equations over $A$ that is solvable in $B$ has a solution in $A$, we say that $A$ is a \emph{pure} subact of $B$. An act is called \emph{absolutely pure} if it is pure in every containing act. This is equivalent to saying that $A$ is injective with respect to all inclusions $K \rightarrow L$ with $L$ a finitely presented and $K$ a finitely generated $S$-act \cite{Normak}. If we assume that $L$ is cyclic in the above statement, we get the concept of \emph{almost pure} acts \cite{Gould87}. For a cardinal $\alpha$, an act is called \emph{$\alpha$-injective} \cite{Gould85} if it is injective with respect to all inclusions $I \rightarrow S$ where $I$ is a right ideal of $S$ the cardinality of which is less than $\alpha$. If $\alpha = \aleph_0$ ($\alpha = 2$) we get the concept of \emph{weakly f-injective} (\emph{weakly p-injective}) acts, respectively \cite{Gould85}.
For any background information or undefined terms about acts over a monoid the reader is referred to \cite{Ahsan, KKM}.
 
Preenvelopes and precovers of modules were introduced by E. Enochs \cite{Enochs} as generalizations of injective envelopes and projective covers, respectively. Analogously, they can be defined for act categories as follows. Let $\mathscr{C}$ be a class of acts over a monoid $S$ and $A$ an $S$-act. A map $\varphi : A \rightarrow C$, with $C \in \mathscr{C}$, is called a $\mathscr{C}$-\emph{preenvelope} of $A$ if for any map $f : A \rightarrow C'$, with $C' \in \mathscr{C}$, there is a map $g: C \rightarrow C'$ such that $g \circ \varphi = f$. A $\mathscr{C}$-\emph{envelope} of $A$ is a $\mathscr{C}$-preenvelope with the added condition that if $C = C'$ and $f = \varphi$ then $g$ must be an automorphism. Hence, a $\mathscr{C}$-envelope of an act, if it exists, is unique up to isomorphism. The class $\mathscr{C}$ is called \emph{(pre)enveloping} if every act admits a $\mathscr{C}$-(pre)envelope. Covers and precovers are defined dually (see, e.g. \cite{BR}). Injective envelopes, defined in \cite{Bert}, are easily seen to be $\mathscr{I}$-envelopes, where $\mathscr{I}$ is the class of injective $S$-acts. It is easy to see that if the class $\mathscr{C}$ contains the injectives then $\mathscr{C}$-preenvelopes are monic. 

The set-theoretic approach to the existence of preenvelopes of acts is adapted from \cite{EandJ} and \cite{RandS}.
\section{Preenvelopes of Acts}
\noindent A sufficient condition for the existence of preenvelopes is given in the following proposition. In what follows, $\abs{X}$ means the cardinality of a set $X$.
\begin{proposition} \label{propn622} Let $A$ be an $S$-act and $\mathscr{C}$ a class of $S$-acts closed under direct products. Suppose that there is an infinite cardinal $\aleph_\alpha$ such that for every $C \in \mathscr{C}$ and every subact $T$ of $C$ with cardinality not exceeding $\abs{A}$ there exists a subact $B$ of $C$ containing $T$ such that $B \in \mathscr{C}$ and $\abs{B} \leq \aleph_\alpha$. Then $A$ has a $\mathscr{C}$-preenvelope.
\begin{proof} Fix a set $X$ with cardinality $\aleph_\alpha$. If we can prove that any map $f: A \rightarrow C$, where $C \in \mathscr{C}$, can be factored through some $B \in \mathscr{C}$ such that $\abs{B} \leq \aleph_\alpha$, then we can make a copy of $B$ with elements from the above set $X$. Hence, the totality of all maps $A \rightarrow B$ in the factorization $A \rightarrow B \rightarrow C$ is in fact a set, and the map $A \rightarrow \prod  B$ is a $\mathscr{C}$-preenvelope. But for the map $f: A \rightarrow C$, we have $T = f(A)$ and $\abs{T} \leq \abs{A}$. Such $T$ is included, by assumption, in a $B$ as above, and we get our desired factorization $A \rightarrow B \rightarrow C$.
\end{proof}
\end{proposition}

Before proceeding to our main result, we need the following lemma, which shows that
between any act and a subact, we can always find a pure one whose cardinality depends only on that of the monoid $S$ and the subact.
\begin{lemma} \label{lem} Let $T$ be a subact of an $S$-act $A$. Then there exist a cardinal $\aleph_\alpha$ depending on the cardinalities of $S$ and $T$, and a pure subact $U$ of $A$ containing $T$ such that $\abs{U} \leq \aleph_\alpha$. 
\begin{proof} Let $\aleph_\alpha = \max\{\aleph _0, \abs{S}, \abs{T}\}$. Consider the set $\Sigma$ of all finite systems of equations of the form 
$$x_ir_i=x_is_i, \thinspace \thinspace y_jr'_j=z_js'_j,\thinspace \thinspace w_kr''_k=t_k; \thinspace \thinspace 1 \leq i \leq  l, \thinspace \thinspace 1 \leq j \leq m, \thinspace \thinspace 1 \leq k \leq n$$
with $r$'s and $s$'s $\in S$, $t$'s $\in T$ and indefinite $x$'s, $y$'s, $z$'s, and $w$'s, that are solvable in $A$. For each variation of $l$, $m$, and $n$ we have a subset $\Sigma_{l,m,n}$ of $\Sigma$ with cardinality $\abs{S}^l .\abs{S}^l. \abs{S}^m .\abs{S}^m .\abs{S}^n .\abs{T}^n = \abs{S}^{2l+2m+n}. \abs{T}^n \leq \aleph_\alpha^{2l+2m+n}. \aleph_\alpha^n =\aleph_\alpha$. Therefore, the cardinality of $\Sigma$ is $\abs{\Sigma} \leq \aleph_0.\aleph_\alpha = \aleph_\alpha$.

Put $U_0 =T$ and consider the set $X_1$ of all systems in $\Sigma$ above with constants from $U_0$ that are solvable in $A$. For each such system, choose one solution. Let $\xbar{X_1}$ be the set of these solutions. Therefore, $\abs{\xbar{X_1}} \leq \abs{\Sigma} \leq \aleph_\alpha$. Now let $U_1 =U_0 \cup \xbar{X_1}S$. Hence $\abs{U_1} \leq \abs{U_0} + \abs{\xbar{X_1}}\abs{S} \leq \aleph_\alpha + \aleph_\alpha . \aleph_\alpha = \aleph_\alpha$.
Similarly, we construct the acts $U_2$, $U_3$, $\cdots$, each of which having cardinality not exceeding $\aleph_{\alpha}$ such that $U_0 \subseteq U_1 \subseteq U_2 \subseteq \cdots \subseteq A$ and that every finite system of equations over $A$ that is solvable in $U_n$ has a solution in $U_{n+1}$. So, putting $U = \bigcup^\infty _{n=0}U_n$, we see that any finite system of equations with constants from $U$ that is solvable in $A$ must have its constants  in some $U_n$. Therefore the system must be solvable in $U_{n+1} \subseteq U$. This proves that $U$ is pure in $A$.

Now for $\abs{U}$. We know that $U = \bigcup^\infty _{n=0}U_n$, and therefore $\abs{U} \leq \abs{U_0} +\abs{U_1} + \cdots \leq \aleph_\alpha + \aleph_\alpha + \cdots = \aleph_0 . \aleph_\alpha = \aleph_\alpha$.
\end{proof}
\end{lemma}

\begin{theorem} \label{main} Let $\mathscr{C}$ be a class of acts closed for pure subacts. Then, $\mathscr{C}$ is preenveloping if and only if it is closed for direct products.
\begin{proof} \emph{Sufficiency:} let $A$ and $C$ be  $S$-acts such that $C \in \mathscr{C}$ and let $T$ be a subact of $C$ with $\abs{T} \leq \abs{A}$. By lemma \ref{lem} and its proof, there is a pure subact $U$ of $C$ containing $T$ with $\abs{U} \leq \max \{\aleph_0, \abs{T}, \abs{S}\} \leq \max \{\aleph_0, \abs{A}, \abs{S}\}$ and, by assumption, $U \in \mathscr{C}$. So, putting $\aleph_\alpha = \max \{\aleph_0, \abs{A}, \abs{S}\}$, we see that $\abs{U} \leq \aleph_\alpha$ and that $U$ is the act $B$ of proposition \ref{propn622}. Hence $A$ has a $\mathscr{C}$-preenvelope by proposition \ref{propn622}. \emph{Necessity:} let $\{ C_i \}$ be a family of acts in $\mathscr{C}$ and let $f: \prod C_i \rightarrow C$ be a $\mathscr{C}$-preenvelope of $\prod C_i$. This means that for every $C_i$ there is a $\lambda _i :C \rightarrow C_i$ such that $\lambda _i f = \pi _i$, where $\pi _i$ is the projection of $\prod C_i$ onto $C_i$ for each $i$. Hence, there is a (unique) $\lambda : C \rightarrow \prod C_i$ such that $\pi _i \lambda = \lambda _i$ for all $i$. But then $\lambda f$ is the identity map of $\prod C_i$, and $\prod C_i$ is a retract, hence a pure subact, of $C$. Therefore, $\prod C_i \in \mathscr{C}$.
\end{proof}
\end{theorem}

\begin{corollary} Each of the classes of absolutely pure, almost pure, weakly f-injective and weakly p-injective acts is preenveloping. Such preenvelopes are all monomorphisms.
\begin{proof} It is easily seen that the above classes are closed for direct products and pure subacts. Therefore, by Theorem \ref{main}, they are all preenveloping classes. Moreover, as each class of them contains the injectives as a subclass, such preenvelopes must then be monic.
\end{proof}
\end{corollary}

\bibliographystyle{amsplain}

\end{document}